\documentclass [12pt]{article}
\usepackage {graphicx}
\usepackage{amsmath,amssymb,newlfont}
\usepackage[T2A]{fontenc}
\usepackage[cp1251]{inputenc}
\usepackage[russian, english]{babel}
\newtheorem{theorem}{\bf Theorem}[section]

\newtheorem{remark}{\bf Remark}[section]

\begin{document}

\vskip 0,5mm
        \begin{center}{\bf The boundary value problem for second order ordinary linear differential
        equations with variable coefficients}\\
\bigskip
        \em A.~Tungatarov, D.K.~Akhmed-Zaki, S.A.~Abdymanapov\\
        al-Farabi Kazakh National University, Almaty, Republic of
        Kazakhstan\\
         e-mail: tun-mat@list.ru\end{center}

\medskip

\noindent {\it 2010MSC}:34С30
\bigskip

\begin{abstract}
      In this work two-point boundary value problem for one class
      of second order ordinary differential equations with
      variable coefficients is solved.\end{abstract}

\vskip 3mm

\section{Introduction}
\label{Introduction}

Let   $0<x_{1}<\infty$, $S[0 ,\, x_{1}]$ is the class of
measurable, essentially bounded functions $f(x)$ in  $[0,\,
x_{1}]$ and  $W_{\infty}^{2}[0,x_{1}]$ is the class of functions
$f(x)$, for which   $\frac{d^2 f}{d x^{2}}\in S[0, \, x_{1}].$
The norm of an element from $S[0,\, x_{1}]$ is defined by the
formulas
$$|f|_{0}=essup_{x\in[0 ,x_{1}]}|f(x)|=\lim\limits_{p\to\infty}\|f\|_{L_{p}[0 ,\, x_{1}]}.$$

We consider the equation

\begin{equation}\label{1}\frac{d^{2}u}{dx^{2}}+a(x)u=f(x)\end{equation}
in interval $[0, \, x_{1}]$, where   $a(x),$ $f(x)\in S[0,\,
x_{1}].$

Two-point boundary value problems for second order ordinary
differential  equations  are classical area of research of the
theory of ordinary differential equations and because of   their
broad application in mechanics, mathematical physics and geometry
(see, for example,\cite{Kam}-\cite{Er}) they are still actively
investigated. However, in mathematical literature the equations of
the form  \eqref {1} with continuous coefficients are studied and
sufficient conditions of resolvability of boundary value problems
for them are received.
 In  author's works  \cite {Tun},
\cite {Cau}  the general solution of equation \eqref {1} is
constructed and  Cauchy problem for it  with initial point $x=0$
is solved. In this work  an explicit form of  general solution of
equation \eqref {1} in class

\begin{equation}\label{2}W_{\infty}^{2}[0, x_{1}]\bigcap C^{1} [0,\, x_{1}],\end{equation}

where \begin{equation}\label{19} x_{1} <
\sqrt{\frac{2}{|a|_{0}}}\end{equation}

is found  and next boundary value problem is  solved.

\textbf{Problem D.} {Find the solution of equation \eqref{1} from
the class   \eqref{2} satisfying the conditions}
\begin{equation}\label{3}u(0)=\alpha,\;\;
u'(x_{1})=\beta,\end{equation}  where $\alpha,  \beta$ {are given
real numbers.}


\vskip 2mm

\section{Construction of the general solutions to equation
\eqref{1}} \label{Tung:sect2}

By integrating two times the equation \eqref{1}, we get

\begin{equation}\label{4}u(x)=(Bu)(x)+g(x)+c_{1} x + c_{2},\end{equation}
where $c_{1}, c_{2}$ are any real numbers,

$$(Bu)(x)=\int\limits_{{0}}^{{x}}\int\limits_{y}^{x_{1}} a(t)u(t)dtdy, \;\;
g(x)=\int\limits_{{0}}^{x}\int\limits_{y}^{x_{1}}f(t)dtdy.$$

Applying the operator  $B$  to both sides of equation \eqref{4} we
have
\begin{equation}\label{5}(Bu)(x)=(B^{2}u)(x)+(Bg)(x)+c_{1}a_{1}(x)+c_{2}b_{1}(x),\end{equation}

where
$$(B^{2}u)(x)=(B(Bu)(x))(x),$$
$$a_{1}(x) = \int\limits_{0}^{x}\int\limits_{y}^{x_{1}}
ta(t)dtdy,\,\,\,\, b_{1}(x) =
\int\limits_{0}^{x}\int\limits_{y}^{x_{1}} a(t)dtdy.$$

From \eqref{4} and \eqref{5} it follows
\begin{equation}\label{6}u(x)=(B^{2}u)(x)+c_{1}(x+a_{1}(x))+c_{2}(1+b_{1}(x))+g(x)+(Bg)(x).\end{equation}

Further we use following functions and operators:
$$a_{k}(x)=\int\limits_{0}^{x}\int\limits_{y}^{x_{1}}
a(t)a_{k-1}(t)dtdy,\,\,
b_{k}(x)=\int\limits_{0}^{x}\int\limits_{y}^{x_{1}}
a(t)b_{k-1}(t)dtdy, \,\, (k=2,3,...),$$
$$(B^{k}u)(x)=(B(B^{k-1}u)(x))(x), \,\, (k=2,3,...).$$

Applying the  operator  $B$ to both sides of  equation \eqref{6}
we get
\begin{equation}\label{7}(Bu)(x)=(B^{3}u)(x)+c_{1}(a_{1}(x)+a_{2}(x))+c_{2}(b_{1}(x)+b_{2}(x))+(Bg)(x)+(B^{2}g)(x).
\end{equation}

From \eqref{4} and \eqref{7} it follows

$$u(x){=}(B^{3}u)(x)+c_{1}(x+a_{1}(x)+a_{2}(x))+c_{2}(1+b_{1}(x)+b_{2}(x))+g(x)+(Bg)(x)+(B^{2}g)(x).$$

Continuing this procedure $n$ times we  obtain the following
integral representation of  solutions of equation \eqref{1}:

\begin{equation}\label{8}u(x)=(B^{n}u)(x)+c_{1}(x+\sum\limits_{k=1}^{n-1}a_{k}(x))+
 c_{2}(1+\sum\limits_{k=1}^{n-1}b_{k}(x))+\sum\limits_{k=0}^{n-1}(B^{k}g)(x),\end{equation}

where $(B^{0}g)(x)=g(x).$

Let  $u(x)\in C[0, x_{1}]$. The following    inequalities are
easily obtained:
\begin{equation}\label{9}|(B^{n}u)(x)|\leq 2|u|_{1}\cdot \frac{|a|_{0}^{n}\cdot
x_{1}^{2n}}{{2}^{n}},\;\;(n=1,2,...),
\end{equation}

\begin{equation}\label{10} |a_{k}(x)| <\ \frac{|a|_{0}^{k}\cdot x_{1}^{2k}}{{2}^{k}}\cdot x,\;\;
 |b_{k}(x)| <\ \frac{2\cdot|a|_{0}^{k}\cdot x_{1}^{2k}}{{2}^{k}},\,\, (k=1,2,...),\end{equation}
where $|f|_{1}=\max_{x\in[0,\, x_{1}]}|f(x)|.$

Passing to the  limit  with $n\to\infty$ in the representation
\eqref{8} and  taking inequalities \eqref{9}, \eqref{19}  into
account
 we get

\begin{equation}\label{12}u(x)= c_{1}I_{1}(x)+ c_{2}I_{2}(x)+F(x),
\end{equation}

where
$$I_{1}(x)=x+\sum\limits_{k=1}^{\infty}a_{k}(x),\;\;I_{2}(x)=1+\sum\limits_{k=1}^{\infty}b_{k}(x),\;\;
F(x)=\sum\limits_{k=0}^{\infty}(B^{k}g)(x).$$

Using the inequalities  \eqref{9}, \eqref{10}, we receive

\begin{equation}\label{13}|I_{1}(x)|\leq\frac{2}{2-|a|_{0}\cdot x_{1}^{2}},\;\;
|I_{2}(x)|< \frac{2+|a|_{0}\cdot x_{1}^{2}}{2-|a|_{0}\cdot
x_{1}^{2}},\;\; |F(x)|\leq|g|_{1}\frac{2+|a|_{0}\cdot
x_{1}^{2}}{2-|a|_{0}\cdot x_{1}^{2}}.\end{equation}

From the form of functions  $I_{1}(x),I_{2}(x)$ and $F(x)$ for
$x\in[0,\, x_{1}]$, where $x_{1}$  satisfies  inequality
\eqref{19}, it follows

\begin{equation}\label{14}\begin{array}{l}\displaystyle
I'_{1}(x)=1+\int\limits_{x}^{x_{1}}a(t)I_{1}(t)dt,\;\; I'_{2}(x)=\int\limits_{x}^{x_{1}}a(t)I_{2}(t)dt, \\\\
\displaystyle
F'(x)=-\int\limits_{x}^{x_{1}}f(t)dt+\int\limits_{x}^{x_{1}}a(t)F(t)dt,\end{array}\end{equation}

\begin{equation}\label{15} I''_{1}(x)=-a(x)I_{1}(x),\;\; I''_{2}(x)=-a(x)I_{2}(x),\end{equation}

\begin{equation}\label{16} F''(x)=f(x)-a(x)F(x).\end{equation}

Considering   \eqref{14} and the form of functions  $I_{1}(x),
I_{2}(x), F(x)$ we get

\begin{equation}\label{17}I_{1}(0)=F(0)=F'(x_{1})=I'_{2} (x_{1})=0, \;\;I_{2}(0)=I'_{1}(x_{1})=1.\end{equation}

From \eqref{15}  and \eqref{16} it follows, that functions
$I_{1}(x), I_{2}(x)$ are  particular solutions from class
\eqref{2} of  homogeneous equations

$$\frac{d^{2}u}{dx^{2}}+a(x)u=0,$$
and  the function $F(x)$ is solution of  non-homogeneous equation
\eqref{1}.

From \eqref{15} and \eqref{17} we see that the Wronskian $W(x)$ is
equal to $-I_{2}(x_{1}).$ Therefore  if $I_{2}(x_{1})\neq 0$ then
the functions $I_{1}(x)$ and $I_{2}(x)$ are linear independent on
$[0, \, x_{1}]$ and the general solution to equation \eqref{1} is
determined by the formula \eqref{12}.

Hence, we proved the following theorem.

\begin{theorem} If $I_{2}(x_{1})\neq 0$ then the function $u(x)$, given by the formula
\eqref{12}, is a general  solution of equation  \eqref{1} from
class \eqref{2}.
\end{theorem}

\section{Solution of boundary value problem}

 For equation   \eqref{1} we consider the problem D.  To solve the problem
 D we use the solution of equation
  \eqref{1}, given by the formula  \eqref{12}. Substituting the function $u(x)$,
 given by formula  \eqref{12}, into boundary conditions  \eqref{3} and taking into account
  \eqref{17} we have

$$c_{2}=\alpha,\;\; c_{1}=\beta.$$
Hence, the solution of problem D has a form
\begin{equation}\label{18}u(x)=\beta I_{1}(x)+\alpha I_{2}(x)+F(x).\end{equation}

Therefore, the following theorem is proved.

\begin{theorem} If $I_{2}(x_{1})\neq 0$ then the problem  D
has solution, which is given by the formula \eqref{18}.
\end{theorem}

\begin{remark} Obviously, that the results of present work  and
take place for  $a(x), f(x)\in C[0,\, x_{1}]$. In this case the
solutions  given by the formulas \eqref{12} and \eqref{18} belong
to the class  $C^{2}[0, \, x_{1}]$. \end{remark}

\end{document}